# Rainfall Forecasting Using Wavelet Transform (Case Study: Rainfall period of West Seram District of Maluku Province-Indonesia)

F. Y. Rumlawang[1, b)] and H. Batkunde[2, a)]

[1,2]*Department of Mathematics, Mathematics and Natural Science Faculty, Pattimura University, Ambon-Inndonesia*
*a) Corresponding author: rumlawang@fmipa.unpatti.ac.id*
*b h.batkunde@fmipa.unpatti.ac.id)*

**Abstract.** West Seram District is an area of Maluku Province that has high rainfall intensity. Trend of erratic weather changes nowadays makes it important to know the information about rainfall intensity. To know this information, we need to know about rainfall period. This will give some contributions to forecast the rainfall later. Rainfall forecasting plays an important role in water resources management. It also can be used to control the unusual event related to rainfall. In this paper, we use Wavelet transform to analyze rainfall data of West Seram District of Maluku Province-Indonesia. The data that we use is data of monthly rainfall from 1991-2016 from Meteorology, Climatology, and Geophysics agency in Ambon-Maluku-Indonesia. We also show some simulation results using Matlab.

## INTRODUCTION

Weather is the state of the atmosphere. This state will continue to fluctuate irregularly and it will difficult to predict the state precisely. The weather information is important to support various activities in many sectors. One of the information can be rainfall forecasting. Rainfall is the amount of rain that falls in a place (within a given area) during a particular period. It includes rain drop, snow, sleet, and dew. The information about rainfall intensity and rainfall period can be used to control many activities such as: agricultural and plantation productions, fisheries, aviation, public services, etc. This is necessary for control arrangements for the above activities in a country, a province or a district.

West Seram District is one of districts in Maluku Province-Indonesia that has high rainfall. It also has high intensity of rain. This district is one of the largest producers of vegetables and rice in Maluku. Local people get their income from this sector. That is why rainfall and rainfall period of this district is important to be known. This information is useful to determine the right planting time. Trends of erratic current weather changes makes this information becomes even more important. It can be used to maintain control, so all the activities of the local people can be accelerated.

To determine rainfall period in an area, we can use some methods in signal processing [1]. One of the most popular methods that can be used to determine rainfall period is Fourier transform. It gives some information about frequency of time series data, but it does not give the time when that frequency happened. This is the shortcoming of Fourier transform. It can be overcome if we use Wavelet transform [2]. The transform can be used to detect time and frequency simultaneously. It is effective to use Wavelet transform to determine rainfall period of an area [3-5]. This wavelet transform be used in this paper to determine rainfall period of West Seram District. The result are expected to help activities of local people.

The theory of wavelet is a concept that relatively new. This concept was introduced by Jean Morlet in 1982 [6]. This French Geophysicist named the concept "wavelet" which come from French "ondelette" that means little wave. Study of wavelet wavelet transform is a new tool for seismic signal analyzing. Then a theoretical physicist Alex Grossman studied inverse formula for the wavelet transform. Both of them collaborated and produced a detailed mathematical study of continuous wavelet transform and its applications. Actually the same concept had already been found in 1950's by Caldeorn, Littlewood, and Paley. The rediscovery of this old concept provides a new method to decompose signals or functions. Its modern applications such as wave propagation [7], data compression [8], signal processing [9], image processing [10], pattern [11], recognition [12], computer graphics[13], etc.

Wavelet analysis is originally introduced in order to improve seismic signal analysis by switching from shortime Fourier analysis to new better algorithms to detect and analyze abrupt changes in signal. The wavelet transform was introduced to circumvent the Heisenberg uncertainty by using adjustable-width basis functions called wavelet that could stretch or compress depending on whether the wavelet is on a low-frequency section or a high-frequency section of the signal being analyzed. A wavelet has a finite duration and zero mean. A function $\psi$ be with

$$\int_{-\infty}^{\infty} \psi(x)\, dx = 0,$$

Is called a wavelet. For every $f, \psi$ defines the continuous wavelet transform

$$W_\psi f(a, b) = \int_{-\infty}^{\infty} f(x) \frac{1}{\sqrt{a}} \overline{\psi\left(\frac{x-b}{a}\right)} dx, \text{ for all } a, b \in \mathbb{R}_+ \times \mathbb{R},$$

where $\psi_{a,b} = \frac{1}{\sqrt{a}} \psi\left(\frac{x-b}{a}\right)$. The function $\psi$ is the so called mother wavelet. It is chosen to be localized at $x = 0$ and some $\omega = \omega_0 > 0$ and/or $\omega = -\omega_+$ [elin Johansen. Litentiate thesis, lulea university of technology department of math 2005]

We also know a concept of short wavelet. Let $\phi$ be a Haar function where

$$\phi(x) = \begin{cases} -1, & -1 < x \leq 0; \\ 0, & 0 < x \leq -1; \\ 1, & \text{else}; \end{cases}$$

transformed into

$$\psi^H(x) = \begin{cases} -\frac{1}{\sqrt{2}}, & -1 < u \leq 0; \\ \frac{1}{\sqrt{2}}, & 0 < u \leq -1; \\ 0, & \text{else}; \end{cases}$$

then $\psi^H$ is a wavelet function. The function $\psi^H$ is called Haar basic wavelet function.

Moreover, wavelet transform is a linier transform which a bit similar to Fourier transform. The fundamental difference between the two is Wavelet allows put time in frequency components which differ from given signal. Wavelet transform is not time-bounded, it can be applied to time series data. Rainfall data of an area is an example of time series data. Time series data is a function of time and among its observations, we can see some relations which called autocorrelations.

## RESULTS

To determine rainfall period in West Seram District, we use daubechies wavelet transform order 1 or known as Haar wavelet. Wavelet with Haar filter is chosen because it has low pass filter and high pass filter that more efficient. In this stage, the signal in form of rainfall data is processed using Haar Wavelet and divided into four levels. The simulation results (using Matlab k2009a) shows that each level produce a period of the signal. Moreover, to reduce some noises in this process, we use thresholding process.in Figure 1 shows a view of initial coefficient using Haar wavelet that divided into four levels using fixed form threshold method.

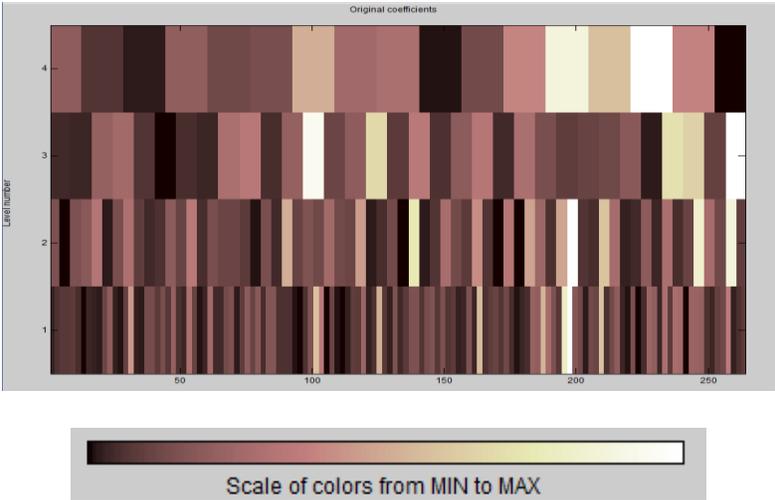

**FIGURE 1.** Original coefficient

The fixed form threshold method is chosen as the best thresholding method in this process because the noises is greatly reduce using this method. Furthermore, we use hard fixed form threshold with non-white noise as its noise structure selection. The following figure shows the view of its de-noise.

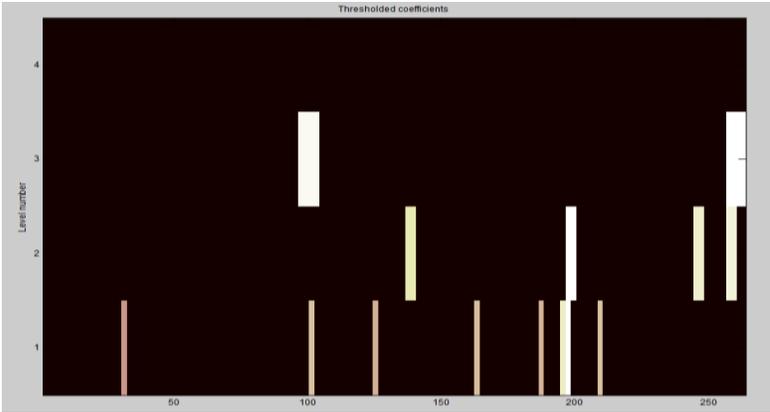

**FIGURE 2.** Thresholded coefficient

Moreover, we determine each level of decomposition so we reduce the noise as we can see in Figure 3.

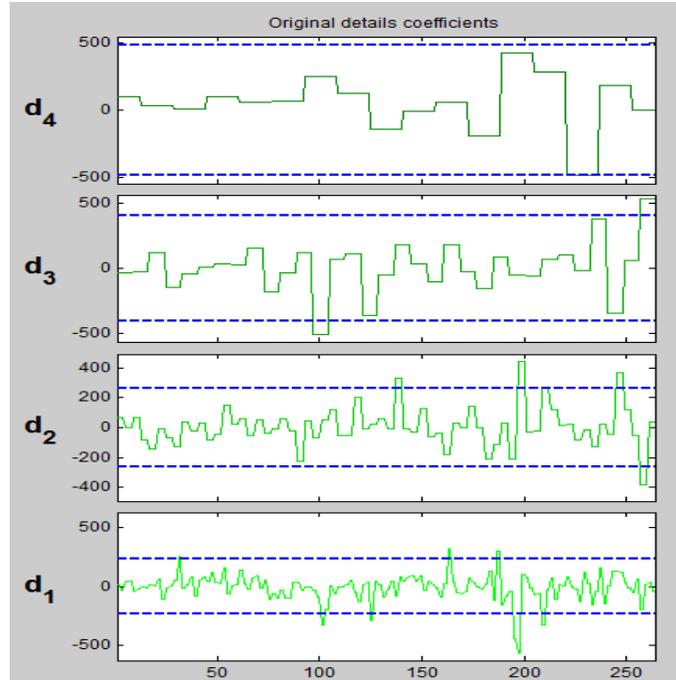

**FIGURE 3.** The thresholding process for each level

Figure 3 shows that the value of these coefficients is above the thresholding value start from first to third level. On the fourth level, the value of the coefficient is below the threshold value. We consider this as a removable noise. Then we chose level three as the highest level.

By using the result of thresholding signal, we can determine the range of time area and range of frequency of the signal. This contains in the following table

**TABLE 1.** Period and frequency of rainfall for West Seram Disrict

| Level | Range | | Median | | Data range of rainfall | Data frequency |
|---|---|---|---|---|---|---|
| | Period (month) | Relative frequency | Period (Month) | Relative frequency | | |
| 1 | 1 – 2 | 0,5 – 1 | 2 | 0,5 | 230 – 950 | 59 |
| 2 | 2 – 4 | 0,25 – 0,5 | 3 | 0,333 | 263 – 950 | 44 |
| 3 | 4 – 8 | 0,125 – 0,25 | 6 | 0,166 | 406 – 950 | 18 |

Based on Table 1, one can see that rainfall period is divided into 3 levels. Rainfall period and frequency is obtained from median value of period range and relative frequency. This median value (rounds up to the nearest integer) will be used as period in each level. Frequency is the number of occurrences of a repeating event per unit of time. We get the relative frequency on Table 1, by formula one over period. On the other hand, period range is obtained from incoming data frequency in time zone on each level of wavelet transform.

Moreover, we can see the result of wavelet transform on Figure 1. It shows rainfall data with respect to frequency and scale. Using different colors, it describes the lowest to highest rainfall. Black color shows that the range has lowest rainfall, and white color shows the highest. From this result, one can see the relation between frequency and scale. Small scale values imply highest frequencies and vice versa. This means wavelet frequencies is changing based on scale value.

Using information on Figure 1 and Table 1, we obtain the rainfall period for level 1 is two months with the lowest rainfall is happened on month 13 to 14 or at January to February 1992. That happened again in month 27 to 28 or at March to April 1993. It also happened in month 35 to 36 or at November to December 1993.

Moreover for level 2, with period three months. One can see that the lowest rainfall happened on month 5 to 7 or at Mei to July 1991. It repeats in month 21 to 23 or September to November 1992. It happened again in month 37 to 39 or January to March 1994, also for month 77 to 79 or May to July 1997 and month 133 to 135 or January to March 2002. We can see that it repeats after sixteen months.

Level 3 has six months period. The lowest rainfall happened in month 41 to 47 (May to November 1994). It repeats in month 56 to 64, 81 to 88, 105 to 112, 129 to 136 and 225 to 231. It repeats after 24 months for this cycle.

By using median value of each rainfall range on Table 1, we obtained rainfall period for one year as shown in the following figure

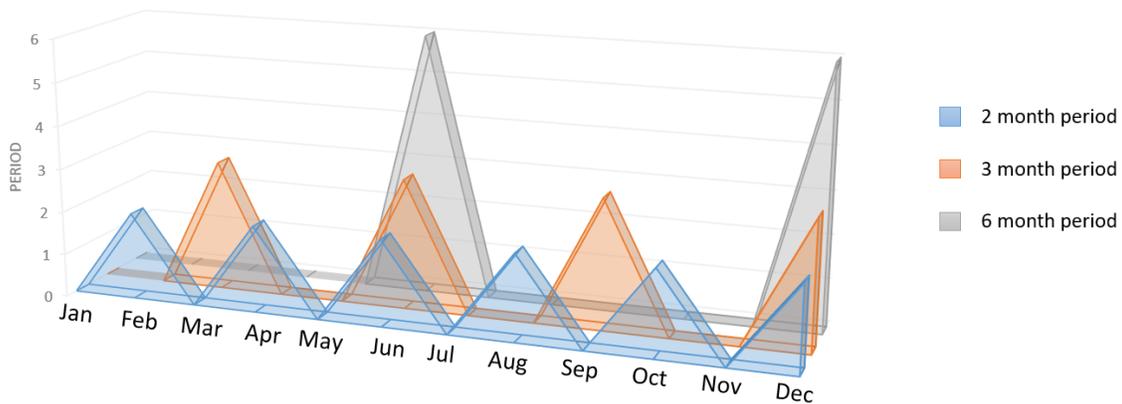

**FIGURE 4.** Rainfall period of West Seram district.

On Fig. 4 one can see that using wavelet transform, we obtain two peak of the rain (bimodial). Those are June and December which is an accumulation of 2, 3, and 6 month periods. The following table gives wavelet transform result for rainfall.

**TABLE 2.** Wavelet transform result

| Period | Months |
|---|---|
| 2 | February, April, June, August, October, December |
| 3 | March, June, September, December |
| 6 | June December |

It can be seen on Table 2 that by using wavelet transform, the rainfall period can be divided into three levels started from the minimum to maximum. Therefore, the transform did the analysis for each frequency based on different resolution. Moreover, by using wavelet transform we have information about timing and frequency of simultaneously rainfall.

## CONCLUSION

Based on the results we can conclude that:
1. By using wavelet transform, we can obtain minimum to maximum period of rainfall for West Seram District. Those periods are two, three, and six month.
2. Wavelet transform shows that the rainfall pattern for West Seram District is equatorial.